\documentstyle[12pt]{article}
\newcommand{\ps}[1]{\sum^{\infty}_{{#1}=1}}
\newcommand{\pzs}[1]{\sum^{\infty}_{{#1}=0}}

\newcommand{\be}{\begin{equation}}
\newcommand{\ee}{\end{equation}}
\newcommand{\bea}{\begin{eqnarray}}
\newcommand{\eea}{\end{eqnarray}}
\newcommand{\kv}[1]{{|#1 \rangle}}
\newcommand{\ket}[1]{{|#1 \rangle}}
\newcommand{\bv}[1]{{\langle #1 |}}

\newcommand{\eq}[1]{\,(\ref{e:#1})\,}
\renewcommand{\v}{${\cal V}ir\hspace{-.03in}_{q,t}\,$}
\newcommand{\vm}{{\cal V}\hspace{-.03in}_{q,t}\ }
\newcommand{\ignore}[1]{}
\newcommand{\n}{\nonumber\\}
\newcommand{\TT}{\tilde{T}}
\renewcommand{\i}{\sqrt{-1}}
\makeatletter

\font\germxxv=eufm10 scaled \magstep5
\@addfontinfo\@xpt{\def\pgerm{\germx}}
\@addfontinfo\@xipt{\def\pgerm{\germxi}}
\@addfontinfo\@xiipt{\def\pgerm{\germxii}}
\@addfontinfo\@xivpt{\def\pgerm{\germxiv}}
\@addfontinfo\@xviipt{\def\pgerm{\germxvii}}
\@addfontinfo\@xxpt{\def\pgerm{\germxx}}
\@addfontinfo\@xxvpt{\def\pgerm{\germxxv}}
\def\germ{\protect\pgerm}
\makeatother
\newcommand\goth[1]{\mbox{\germ{#1}}}
\begin{document}
%
\renewcommand{\thefootnote}{\fnsymbol{footnote}}
\font\csc=cmcsc10 scaled\magstep1
{\baselineskip=14pt
 \rightline{
 \vbox{
       \hbox{YITP/U-95-30}
       \hbox{DPSU-95-5}
       \hbox{UT-715}
       \hbox{July 1995}
}}}

\vskip 5mm
\begin{center}
{\large\bf
A Quantum Deformation of the Virasoro Algebra\\ and\\ \vskip1.6mm
the Macdonald Symmetric Functions}

\vspace{15mm}

{\csc Jun'ichi SHIRAISHI}$^1$,
{\csc Harunobu KUBO}\footnote{JSPS fellow}$^2$,
{\csc Hidetoshi AWATA}$^{*3}$ \\ and \\
{\csc Satoru ODAKE}$^4$
{\baselineskip=15pt
\it\vskip.25in
  \setcounter{footnote}
{0}\renewcommand{\thefootnote}{\arabic{footnote}}\footnote{
      e-mail address : shiraish@momo.issp.u-tokyo.ac.jp}
Institute for Solid State Physics, \\
University of Tokyo, Tokyo 106, Japan \\
\vskip.1in
  \footnote{
      e-mail address : kubo@danjuro.phys.s.u-tokyo.ac.jp}
Department of Physics, Faculty of Science \\
  University of Tokyo, Tokyo 113, Japan \\
\vskip.1in
  \footnote{
      e-mail address : awata@yisun1.yukawa.kyoto-u.ac.jp}
Uji Research Center, Yukawa Institute for Theoretical Physics \\
  Kyoto University, Uji 611, Japan
\vskip.1in
  \footnote{
      e-mail address : odake@yukawa.kyoto-u.ac.jp}
Department of Physics, Faculty of Science \\
  Shinshu University, Matsumoto 390, Japan\\
}
\end{center}

\vspace{10mm}
\begin{abstract}

A quantum deformation of the Virasoro algebra is defined.
The Kac determinants at arbitrary levels are conjectured.
We construct a  bosonic realization of the quantum
deformed Virasoro algebra. Singular vectors
are expressed by the Macdonald symmetric functions.
This is proved by constructing screening currents acting on the bosonic Fock
space.

\end{abstract}

\vspace{7mm}
q-alg/9507034

\vfill
\newpage
%
%
\section{Introduction}
\label{intro}
The hidden correspondence between
the Calogero-Sutherland model (CSM)\cite{rCS}  and the
conformal field theory (CFT)\cite{rBPZ}
is very interesting and seems
even a mysterious thing\cite{rAMOS}\cite{rMY1}.
This correspondence is studied by using the theory of the Jack
symmetric polynomials\cite{rS} and the Feigin-Fuchs
bosonization\cite{rFF} of the Virasoro algebra.
In CSM every excited state is  written in terms of
the Jack symmetric polynomial, and the singular vectors
of CFT are also given by the Jack symmetric
functions with rectangular diagrams.
The Macdonald symmetric functions\cite{rMac} also can be written by
using $q$-analogues of the screening currents
\cite{rAMOS} \cite{rAOS}.
This is a strong evidence that the correspondence
of CSM and CFT still survive after the $q$-deformation.
Therefore an answer to the following question should be desired:
 what algebraic structure dose emerge
from a quantum deformation of
 CSM obtained by replacing the Jack polynomials with
the Macdonald polynomials?
A quite interesting answer to this question does exist,
and we can regard the new algebra as an deformation of the Virasoro
algebra.
In this letter we study
this  quantum deformation of the Virasoro algebra.
Highest weight modules for this quantum deformed
Virasoro ($q$-Virasoro) algebra  are defined.
In order to have reducibility conditions,
we calculate the Kac determinants for
 some lower levels and obtained a conjectural formula
 for general level $N$.
We construct a free boson realization of the $q$-Virasoro algebra.
This bosonization can be regarded as
a $q$-deformation of the Feigin-Fuchs representation of the Virasoro algebra.
We show that, after bosonization, singular vectors are written in terms of the
Macdonald
symmetric functions.
This is analogous to the case of the ordinary
Virasoro algebra whose singular
 vectors are given by  the Jack symmetric
functions.
Surprisingly, our $q$-Virasoro algebra has a connection with the
quantum affine algebra in  the following sense.
Quantum deformation of the Virasoro and ${\cal W}$-algebras
associated with the quantum affine algebras at critical level
were successfully constructed by
Frenkel and Reshetikhin\cite{rFR}.
One can show that our $q$-Virasoro algebra tends to theirs
in a certain limit (see Discussions).
Frenkel and Reshetikhin pointed out their bosonized
$q$-Virasoro and ${\cal W}$-currents resemble
the dressed vacuum form (DVF)
in the algebraic Bethe ansatz\cite{rKS}.
We can notice that the bosonization formula for our $q$-Virasoro algebra
has DVF structure, too (see section 4).
We present a heuristic derivation of our bosonized
$q$-Virasoro current in Appendix.

Applying the collective field theory
\cite{rJS}\cite{rAJL}
to the Calogero-Sutherland model,
 it can be shown  that the Virasoro generators are related to
the Calogero-Sutherland Hamiltonian  $\widehat{\cal {H}}_{\beta}$ as follows
\cite{rAMOS},
\be
\widehat{\cal {H}}_{\beta}= \beta\ps{n} a_{-n}L_n + (\beta(N+1-2 a_0)-1)
                       \widehat{\cal{P}},
\ee
where $\beta$ is the coupling constant of the Calogero-Sutherland model
and $N$
is the number of the particles, $\widehat{\cal{P}}$ is the momentum of CSM,
$L_n$'s are the generators of the Virasoro algebra
and $a_{-n}$'s are the creation operators of the free boson field
 and $a_0$ is its zero mode.
This splitting $\widehat{\cal {H}}_{\beta}$ into $a_{-n}$ and
$L_n$  means that the singular vectors of the  Virasoro algebra
are given by the Jack symmetric functions.
It is natural to expect that we have an  analogous
{\it split expression}
for the bosonized Macdonald operator $\widehat{D}_{q,t}$ \cite{rAMOS}
\be
\widehat{D}_{q,t} = \ps{n} \psi_{-n}T_{n}+\cdots,
\ee
where $T_{n}$'s should be understood
as the generators of the $q$-Virasoro algebra and
$\psi_{-n}$'s
do not contain positive modes of the free boson.
In this article, we will show that this splitting certainly can be done,
and  $T_{n}$'s  generate  an algebra
having the Macdonald symmetric functions as its singular vectors.
We can obtain screening currents $S_{\pm}(z)$ for our $q$-Virasoro algebra.
 Using this free boson realization scheme, we discuss
Feigin-Fuchs like construction of the highest weight
modules of the $q$-Virasoro algebra and study
the singular vectors.
\par
This letter is organized as follows.
In section \ref{definition}, we define the quantum deformation
of the Virasoro algebra.
In section \ref{highest}, the Kac determinant of the
$q$-Virasoro algebra is calculated  at level 1 and level 2.
In section \ref{bosonization}, the free boson realization of
the $q$-Virasoro algebra and  the splitting of
the Macdonald operator are  given. Screening currents are
constructed
in section \ref{screening}, and the singular vectors are written
down by the screening currents. We show the singular vectors are
the Macdonald symmetric functions with rectangular diagrams.
Section \ref{discussion} is devoted to discussions. In
Appendix we show our heuristic derivation of the bosonized
$q$-Virasoro current.
%
%
\section{Definition of the $q$-Virasoro algebra \v }
\label{definition}
Let $q,t$ be two generic complex parameters.
For simplicity  we will consider the case $|p| < 1$,
 where we have set $p=qt^{-1}$ and will use this
notation frequently. We also write $t=q^{\beta}$ by
introducing a complex parameter $\beta$,
which plays the role of the coupling constant
of the CSM in the limit of $q \to 1$.
The $q$-Virasoro algebra \v is an associative algebra generated by
$\{T_n|n\in \bf{Z}\}$ with the following relations
\be
[T_n \, , \, T_m]=-\ps{l}f_l\left(T_{n-l}T_{m+l}-T_{m-l}T_{n+l}\right)
-\frac{(1-q)(1-t^{-1})}{1-p}(p^{n}-p^{-n})\delta_{m+n,0},
\label{e:a1}
\ee
where the coefficients
$f_{l}$'s are given by the following generating function $f(z)$
\be
f(z)=\pzs{l}f_l z^l
=\exp \left\{\ps{n}\frac{1}{n}\frac{(1-q^n)(1-t^{-n})}{1+p^n}
 z^n \right\}.
\label{e:a1.1}
\ee
Introducing \v current
$T(z)=\sum_{n\in\bf{Z}}T_n z^{-n}$, the defining relation \eq{a1} can be
written as follows
\be
f(w/z)T(z)T(w)-T(w)T(z)f(z/w)
= -\frac{(1-q)(1-t^{-1})}{1-p}\left[
       \delta \Bigl(\frac{pw}{z}\Bigr)-
       \delta \Bigl(\frac{p^{-1}w}{z}\Bigr)\right],
\label{e:a1.2}
\ee
where
$\delta(x)=\sum_{n \in {\bf Z}}x^n$.
 Note that the defining relation \eq{a1} is invariant
under the transformation
\be
T_n \to -T_n,
\label{e:a1.3}
\ee
and also invariant under
\be
(q,t)\to (q^{-1},t^{-1}).
\label{e:aa1}
\ee
\par
Here we give the relation between
our quantum deformed Virasoro current $T(z)$
and the ordinary one $L(z)$.
Let us study the limit $q\to 1$ by
parameterizing  $q=e^{h}$. Suppose that $T(z)$ has the expansion
\be
T(z)=2+\beta  \left( z^2 L(z)+\frac{(1-\beta)^2}{4\beta} \right)h^2
           + T^{(2)}(z)h^4 +\cdots.
\label{pe:a7}
\ee
This expansion is consistent with the invariance under transformation \eq{aa1}.
The defining relation \eq{a1} gives us the well known relations for
$ L(z)=\sum_{n \in {\bf Z}}L_{n}z^{-n-2}$, namely
\be
[L_n , L_m]=(n-m)L_{n+m}+\frac{c(n^3-n)}{12}\delta_{n+m,0},
\label{e:a8}
\ee
where
\be
c=1-\frac{6(1-\beta)^2}{\beta}.
\label{e:a9}
\ee
This relation between  the central charge of
the Virasoro algebra and the coupling constant
of CSM is discussed
by studying the Virasoro constraints of generalized
matrix models\cite{rAMOS}.
\section{Highest weight modules of  \v}
\label{highest}
Let us define the Verma module of \v.
Let  $\kv{\lambda}$ be the highest weight vector
having the properties
\bea
& & T_0 \kv{\lambda}=\lambda\kv{\lambda},
\label{e:a14.1}
\\
& & T_n \kv{\lambda}=0   \quad  \mbox{for}  \quad   n\ge 1.
\label{e:a14.2}
\eea
The Verma module $M(\lambda)$ is defined by
\be
M(\lambda)=\vm\kv{\lambda}.
\label{e:a14.3}
\ee
The irreducible highest module $V(\lambda)$ is obtained from
$M(\lambda)$ by removing all singular vectors and their descendants.
Right modules are defined in a similar way from the
lowest weight vector
$\bv{\lambda}$
s.t.
$\bv{\lambda}T_0=\lambda\bv{\lambda},
 \ \bv{\lambda}T_n=0\ n \leq -1$.
A unique invariant paring is defined by setting
$\bv{\lambda}\lambda \rangle = 1$.
The Verma module of the ordinary Virasoro algebra may have the singular
vectors.
The Verma module $M(\lambda)$
may also have the
singular vectors in the same way.
Let us introduce the (outer) grading  operator $d$ which satisfies
$[d , T_n]= n T_n$. Set $d\ket{\lambda}=0$. We call a vector
$\ket{v} \in M(\lambda)$ of level $n$ if $d\ket{v}=-n\ket{v}$.

Whether there exist the singular vectors or not is
checked by calculating the Kac determinant.
Here, we give some explicit forms of $f_{n}$
which we will use for the calculations
\bea
& & f_{1}=\frac{(1-q)(1-t^{-1})}{1+p},
\label{e:a11.2}
\\
& &
f_{2}=\frac{(1-q^{2})(1-t^{-2})}{2(1+p^{2})}+\frac{(1-q)^2(1-t^{-1})^2}{2(1+p)^2}.
\label{e:a11.3}
\eea
\par
At level 1, the Kac determinant is the $1\times 1$
matrix as follows
\be
\bv{\lambda}T_{1} T_{-1}\kv{\lambda}
=\frac{(1-q)(1-t)}{q + t}(\lambda^2 - (p^{1/2}+p^{-1/2})^2 ).
\label{e:a15}
\ee
Therefore, there exist a singular vector at level 1 iff
$\lambda=\pm \left(p^{1/2}+p^{-1/2}\right)$,
 since  $q$ and $t$ are generic.
The signs $\pm$ in the RHS are due to the symmetry \eq{a1.3}.
\par
At level 2, the Kac determinant is
\bea
& &\left|
         \begin{array}{clcr}
              \bv{\lambda}T_{1}T_{1}T_{-1}T_{-1}\kv{\lambda} &
                 \bv{\lambda}T_{1}T_{1}T_{-2}\kv{\lambda} \\
              \bv{\lambda}T_{2}T_{-1}T_{-1}\kv{\lambda} &
                   \bv{\lambda}T_{2}T_{-2}\kv{\lambda} \\
         \end{array}
 \right|
\nonumber \\
&= &\frac{(1-q^2)(1-q)^2 q^{-4}(1-t^2)(1-t)^2 t ^{-4}}{(q +t)^2 (q^2 +t^2)}
\nonumber \\
& &    \times (\lambda^2 qt-(q+t)^2 )( \lambda^2 q^2t-(q^2+t)^2 )
           ( \lambda^2 qt^2 - (q+t^2)^2 ).
\label{e:a18}
\eea
The vanishing conditions of the Kac determinant are
\bea
& \mbox{(i)}&\lambda =\pm \left(p^{1/2} + p^{-1/2}\right),
\label{e:a19.1}  \\
& \mbox{(ii)}&\lambda  =\pm \left(p^{1/2}q^{1/2} + p^{-1/2}q^{-1/2}\right),
\label{e:a19.2}  \\
& \mbox{(iii)}& \lambda =\pm \left(p^{1/2}t^{-1/2} + p^{-1/2}t^{1/2}\right).
\label{e:a19.3}
\eea
In the case (i), there is a singular vector at level 1.
In the cases (ii) and (iii), we have a singular  vector at level 2.
The singular vector for the case (ii) is
\be
 \frac{qt^{-1/2}(q+t)}{(1-q)^2(1+q)}T_{-1}T_{-1}\kv{\lambda}
\mp T_{-2}\kv{\lambda},
\label{e:a20.1}
\ee
and for (iii) is
\be
 \frac{q^{-1/2}t(q+t)}{(1-t)^2(1+t)}T_{-1}T_{-1}\kv{\lambda}
\mp  T_{-2}\kv{\lambda}.
\label{e:a21.1}
\ee
In discussion, we will state a conjecture of the Kac determinant
 for arbitrary level $N$.
%
%
\section{Free boson realization of \v }
\label{bosonization}
In this section we construct  a free boson realization of  \v .
One may expect that a free  boson realization of  \v
enables us to investigate new aspects of the singular vectors.
Let us introduce bosonic oscillators $a_n\ n\in {\bf Z}$ and $Q$ with the
commutation
relations
\be
[a_{n} , a_{m}]=n\frac{1-q^{|n|}}{1-t^{|n|}}\delta_{n+m,0},
\ee
\be
[a_{n},Q]=\frac{1}{\beta}\delta_{n,0}.
\ee
Let us define the bosonic Fock space.
Introduce the vacuum state $\kv{0}$ which satisfy the following conditions
\be
a_{n}\kv{0}=0 \quad\mbox{for}\quad n\ge 0.
\label{e:b3.1}
\ee
Furthermore, add the zero-mode momentum for $r,s \in  {\bf Z}$
\be
\kv{r,s}=e^{\alpha_{r,s}Q}\kv{0},
\label{e:b3.3}
\ee
where
\be
\alpha_{r,s}=\frac{1}{2}(1+r)\beta-\frac{1}{2}(1+s).
\label{e:b3.4}
\ee
We define the Fock space ${\cal F}_{r,s}$ by
${\bf C}[a_{-1},a_{-2},\cdots]\kv{r,s}$.
By studying the operator product expansion,
we can check that \v current $T(z)$ is bosonized as follows:
\bea
T(z)
 & = & p^{1/2}\exp\left\{-\ps{n}\frac{1-t^n}{1+p^n}
\frac{a_{-n}}{n}z^{n}t^{-n}p^{-n/2}\right\}
    \exp\left\{-\ps{n}(1-t^n)\frac{a_{n}}{n}z^{-n}p^{n/2}\right\}
    q^{\beta a_{0}} \nonumber  \\
  & &\!\!+\ p^{-1/2}\exp\left\{\ps{n}\frac{1-t^n}{1+p^n}
\frac{a_{-n}}{n}z^{n}t^{-n}p^{n/2}\right\}
    \exp\left\{\ps{n}(1-t^n)\frac{a_{n}}{n}z^{-n}p^{-n/2}\right\}
    q^{-\beta a_{0}}.
\label{e:b1}
\eea
We can observe that this formula has strong resemblance to the
dressed vacuum form in the algebraic Bethe ansatz.
This profound $q$-Virasoro-DVF correspondence was
discovered by Frenkel and Reshetikhin \cite{rFR}.
In Appendix we state our heuristic derivation of the
formula \eq{b1} which does not rely on the defining relation \eq{a1}.

We have  the following highest weight conditions
\bea
& &T_0\kv{r,s}=\left(p^{1/2}q^{\alpha_{r,s}} +
 p^{-1/2}q^{-\alpha_{r,s}}\right)\kv{r,s},
\nonumber \\
& &
T_n\kv{r,s}=0 \quad\mbox{for}\quad  n\geq 1.
\label{e:b5.1}
\eea
So we obtained the embedding $V(\lambda_{r,s}) \to {\cal F}_{r,s}$,
 where
$\lambda_{r,s}
=\left(p^{1/2}q^{\alpha_{r,s}}+p^{-1/2}q^{-\alpha_{r,s}}\right)$.
\par
Next we will show how the Macdonald operator $D_{q,t}$ is factorized
in terms of the \v-current $T(z)$.
The operator $D_{q,t}$ is defined by
\be
  D_{q,t}=
  \sum_{i=1}^N\prod_{\scriptstyle j=1 \atop\scriptstyle j\neq i}^N
  \frac{tx_i-x_j}{x_i-x_j}T_{q,x_i}\:,
\ee
where $T_{\xi , x_i}$ is the shift operator defined by
$T_{\xi , x_i}g(x_1,\cdots, x_i, \cdots,x_r)=g(x_1,\cdots, \xi x_i,
\cdots,x_r)$.
$D_{q,t}$ acts on the ring of the symmetric polynomials with
 $N$-variables $\Lambda_N$.
We have the projection
$\pi_N: {\bf C}[a_{-1},a_{-2},\cdots]\ket{0} \to
\Lambda_N$ given by the mapping
 $\pi_N: a_{-n_1}a_{-n_2}\cdots\ket{0}\mapsto
p_{n_1} p_{n_2}\cdots
$ where
$p_n= \sum_{i=1}^{N}x_i^n$.
We have the bosonized Macdonald operator $\widehat{D}_{q,t}$
which satisfies $ D_{q,t}\circ \pi_N=\pi_N \circ \widehat{D}_{q,t}$
 given by the following formula\cite{rAMOS}
\be
\widehat{D}_{q,t}=\frac{t^N}{t-1}\oint\frac{dz}{2\pi \i z}
\exp\left\{\ps{n}\frac{1-t^{-n}}{n}a_{-n}z^n\right\}
\exp\left\{-\ps{n}\frac{1-t^{n}}{n}a_{n}z^{-n}\right\}-\frac{1}{t-1}.
\label{e:b6}
\ee
One may find that this bosonized
operator $\widehat{D}_{q,t}$ can be decomposed by using
the $q$-Virasoro current $T(z)$ and the  operator $\psi (z)$ defined by
\be
\psi(z)=\pzs{n}\psi_{-n}z^n=p^{-1/2}\exp\left\{-\ps{n}\frac{1-t^n}{1+p^n}
\frac{a_{-n}}{n}z^np^{n/2}t^{-n}\right\}
                      q^{-\beta a_0 },
\label{e:b7.2}
\ee
as follows,
\bea
\widehat{D}_{q,t}&=&\frac{t^N}{t-1}\left[\oint\frac{dz}{2\pi\i}\frac{1}{z}
\psi(z)T(z) -p^{-1}q^{-2\beta a_0}\right] -\frac{1}{t-1}
\nonumber \\
&= &\frac{t^N}{t-1}\left[\pzs{n}\psi_{-n}T_n
             -p^{-1}q^{-2\beta a_0}\right]-\frac{1}{t-1}.
\label{e:b7.3}
\eea
%
%
\section{Screening charges and singular vectors}
\label{screening}

We study the singular vectors of \v in the bosonic Fock space.
Define  the  screening currents $S_{\pm}(z)$ as follows:
\bea
  S_+(z)
  &\!\!=\!\!&
  \exp\left\{\ps{n}\frac{1-t^n}{1-q^n}\frac{a_{-n}}{n}z^{n}\right\}
  \exp\left\{-\ps{n}(1+p^n)\frac{1-t^n}{1-q^n}\frac{a_{n}}{n}z^{-n}
  \right\} e^{\beta Q}z^{2\beta a_0},
  \label{e:c1.1}\\
  S_-(z)
  &\!\!=\!\!&
  \exp\left\{-\ps{n}\frac{a_{-n}}{n}z^{n}\right\}
  \exp\left\{\ps{n}(1+p^n)\frac{a_{n}}{n}z^{-n}p^{-n}\right\}
  e^{- Q}z^{-2 a_0}.
  \label{e:c1.2}
\eea
The commutation relation between $T_n$ and
the screening currents are
\bea
  \Bigl[T_n,S_+(w)\Bigr]
  &\!\!=\!\!&
  -(1-q)(1-t^{-1})\frac{d_q}{d_q w}
  \left((p^{-\frac12}w)^{n+1}A_+(w)\right), \label{e:c2.1}\\
  \Bigl[T_n,S_-(w)\Bigr]
  &\!\!=\!\!&
  -(1-q^{-1})(1-t)\frac{d_t}{d_t w}
  \left((p^{\frac12}w)^{n+1}A_-(w)\right),
  \label{e:c2.2}
\eea
where
\bea
  A_+(z)
  &\!\!=\!\!&
  \exp\left\{\ps{n}\frac{1+t^n}{1+p^n}\frac{1-t^n}{1-q^n}
  \frac{a_{-n}}{n} z^nt^{-n}\right\}
  \exp\left\{-\ps{n}(1+t^n)\frac{1-t^n}{1-q^n}\frac{a_{n}}{n}
  z^{-n}p^n\right\} \n
  &&\times
  e^{\beta Q }z^{2\beta a_0}t^{-a_0},
  \label{e:c4} \\
  A_-(z)
  &\!\!=\!\!&
  \exp\left\{-\ps{n} \frac{1+q^n}{1+p^n}\frac{a_{-n}}{n}
  z^nt^{-n}\right\}
  \exp\left\{\ps{n} (1+q^n)\frac{a_{n}}{n}
  z^{-n}p^{-n}\right\}e^{-Q}z^{-2 a_0}t^{a_0},
  \label{e:c5}
\eea
and the difference operator with one parameter
is defined by
\be
\frac{d_\xi}{d_\xi z}g(z)=\frac{g(z)-g(\xi z)}{(1-\xi)z}.
\label{e:c6}
\ee
We can construct the BRS charges in the same way discussed in the
paper\cite{rDFF}.
On the Fock space ${\cal F}_{-r,s}$, we have the
 single-valued BRS current
\be
J_{+}^{(r)}(w)=S_{+}(w)\oint\prod^{r}_{i=2}
\frac{dw_i}{2\pi \i}S_{+}(w_i),
\label{e:c6.1}
\ee
where the integration cycle is chosen to be Felder's one.
Even though the commutator $[T_n,S_+(w)]$ is a ``total difference'',
we can check that $[T_n, \oint dw  J^{(r)}_{+}(w)]=0$ on the Fock space
${\cal F}_{-r,s}$.
 This is because on the  Fock space ${\cal F}_{-r,s}$ the expansion
$[T_n, J^{(r)}_{+}(w)]=\sum_{m \in {\bf Z}} {\cal O}_m w^{-m}$ is well
defined and ${\cal O}_1 = 0$ since the commutator is a total difference.
\par
We also have another BRS current
\be
J_{-}^{(s)}(w)=S_{-}(w)\oint\prod^{s}_{i=2}
\frac{dw_i}{2\pi \i}S_{-}(w_i),
\label{e:c6.2}
\ee
which is single-valued on the Fock space ${\cal F}_{r,-s}$.

The singular vector in ${\cal F}_{r,s}$ is written by using
BRS current as
\be
\kv{\chi_{r,s}}=\oint \frac{dw}{2\pi \i}J^{(r)}_{+}(w)\kv{-r,s}
=\oint \frac{dw}{2\pi \i} J^{(s)}_{-}(w)\kv{r,-s},
\label{e:c8}
\ee
up to normalization.
Using equation \eq{b7.3}, we obtain
\be
\widehat{D}_{q,t}\kv{\chi_{r,s}}=
\left( \sum^{r}_{i=1}t^{N-i}q^s +\sum^{N}_{i=r+1}t^{N-i}\right)
\kv{\chi_{r,s}}.
\label{e:c10}
\ee
This means that the singular vector $\kv{\chi_{r,s}}$  is the eigen vector of
the
Macdonald operator $\widehat{D}_{q,t}$. Namely $\kv{\chi_{r,s}}$
 is proportional to the Macdonald symmetric function
$P_{\lambda}(q,t)$ with the
rectangular Young diagram $\lambda=(s^r)$.
\par
It may be helpful to give another proof
that $\kv{\chi_{r,s}}$
is the bosonized Macdonald symmetric function
without using $q$-Virasoro algebra .
We calculate $\widehat{D}_{q,t}\kv{\chi_{r,s}}$
in the following way.
We have
\be
\kv{\chi_{r,s}} =
\oint\prod^{r}_{i=1}\frac{dw_i}{2\pi \i}
F^{(+)}(w_1, w_2, \cdots , w_r)
: S_{+}(w_1)S_{+}(w_2)\cdots S_{+}(w_r)e^{\alpha_{-r,s}Q}:\kv{0},
\label{e:c10.1}
\ee
where
\be
F^{(+)}(w_1,\cdots, w_r)=
\left[
      \prod^{r}_{j=1}w_j^{2\alpha_{-r,s}+2(r-j)\beta}
\right]
\left[
      \prod_{1 \leq i < j \leq r}
      \frac{\left(w_j/w_i;q\right)_{\infty}
            \left(pw_j/w_i;q\right)_{\infty}}
           {\left(tw_j/w_i;q\right)_{\infty}
            \left(qw_j/w_i;q\right)_{\infty}}
\right],
\label{e:c10.2}
\ee
here we used the standard notation $\left(z;q\right)_{\infty}
\equiv \prod_{k=0}^{\infty}\left(1-q^{k}z\right)$.
An explicit residue calculation gives us the following
\bea
& &t^{-N+1}\left(\widehat{D}_{q,t}\  -
\frac{1-t^{N-r}}{1-t}\right)\kv{\chi_{r,s}}
\nonumber \\
& &
=\oint\prod^{r}_{i=1}\frac{dw_i}{2\pi \i}
\sum^{r}_{i=1}q^{-1}\prod_{j=1}^{i-1}t^{-1}\frac{1-p^{-1}w_i/w_j}{1-q^{-1}w_i/w_j}
\cdot\prod_{j=i+1}^{s}\frac{1-pw_j/w_i}{1-qw_j/w_i}\cdot
\left(T_{q^{-1},w_i}
F^{(+)}(w_1,\cdots, w_r)\right)
\nonumber \\
& &\quad  \times :S_{+}(w_1)S_{+}(w_2)\cdots S_{+}(w_r)
e^{\alpha_{-r,s}Q}:\kv{0}
\nonumber \\
& &
=\oint\prod^{r}_{i=1}\frac{dw_i}{2\pi \i}
t^{1-r}q^s\sum^{r}_{i=1}\prod_{i \neq j}\frac{1-tw_j/w_i}{1-w_j/w_i}
F^{(+)}(w_1,\cdots, w_r) \n
&& \quad  \times
:S_{+}(w_1)S_{+}(w_2)\cdots S_{+}(w_r)e^{\alpha_{-r,s}Q}:\kv{0}
\nonumber \\
& &
=\sum^{r}_{i=1}q^{s-(i-1)\beta}\kv{\chi_{r,s}}.
\label{e:c10.4}
\eea
In the calculation we used the following equations
\bea
& &T_{q^{-1},w_i}F^{(+)}(w_1,\cdots, w_r)
\nonumber \\
& &=
\left[
\prod_{j=1}^{i-1}
\frac{\left(1-q^{-1}w_i/w_j\right)
            \left(1-t^{-1}w_i/w_j\right)}
           {\left(1-p^{-1}w_i/w_j\right)
            \left(1-w_i/w_j\right)}
\right]
\left[
\prod_{j=i+1}^{r}
\frac{\left(1-tw_j/w_i\right)
            \left(1-qw_j/w_i\right)}
           {\left(1-w_j/w_i\right)
            \left(1-pw_j/w_i\right)}
\right]
\nonumber \\
& &
\quad \times q^{-2\alpha_{-r,s}}t^{-2(r-i)}F^{(+)}(w_1,\cdots, w_r),
\label{e:nc10.5}
\eea
and
\be
\sum^{r}_{i=1}\prod_{i \neq j}\frac{1-tw_j/w_i}{1-w_j/w_i}=
\frac{1-t^r}{1-t}.
\label{e:c10.6}
\ee
Then we get $
\widehat{D}_{q,t}\kv{\chi_{r,s}}=
\left( \sum^{r}_{i=1}t^{N-i}q^s +\sum^{N}_{i=r+1}t^{N-i}\right)
\kv{\chi_{r,s}}$.
A similar calculation can be done for the singular
vectors composed of  the screening currents $S_{-}(w)$'s.
%
%
\section{Discussions}
\label{discussion}
In this section we look through some of our further results and
unclarified problems.
\begin{itemize}
\item
To calculate the Kac determinant becomes difficult task
when $N$ increases.
We have calculated up to level 4, and write down the
conjectural form at level $N$  as follows
\[
  \det{}_N
  =
  \det\Bigl(\langle i\ket{j}\Bigr)_{1\leq i,j\leq p(N)}
  \!\!=\!\!
  \prod_{\scriptstyle r,s\geq 1 \atop \scriptstyle rs\leq N}
  \Bigl(\lambda^2-\lambda_{r,s}^2\Bigr)^{p(N-rs)}
  \left(\frac{(1-q^r)(1-t^r)}{q^r+t^r}
  \right)^{p(N-rs)},
\]
where the basis at level $N$ is defined
$\ket{1}=T_{-N}\ket{\lambda}$,
$\ket{2}=T_{-N+1}T_{-1}\ket{\lambda}$,$\cdots,\ket{p(N)}=T_{-1}^N\ket{\lambda}$, and $p(N)$ is the
number of the partition of $N$.
We remark that the $\lambda$ dependence has  essentially
the same structure as the case of the usual Virasoro algebra.
Therefore, if $q$ and $t$ are generic,  the character of the quantum Virasoro
algebra \v,
which counts the degeneracy at each level,
exactly coincides with that of the usual Virasoro algebra.
The $\lambda$-independent factor
in the RHS will play an important role when we  study the case
that $q$ is a root of unity.
\item
We can study the  limit of  $\beta \to  0$ and $q$ be  fixed.
In this case $T_n$'s become commutative
and  the defining relations of the $q$-Virasoro algebra\eq{a1}
reduces to the eq. (9.4) of Frenkel-Reshetikhin \cite{rFR}
 if we replace the commutation bracket $\frac{1}{\beta}[\ \ ,\ \ ]$
by the Poisson bracket $\{\ \ ,\ \ \}$.
The limit of $q \to 0$ and $t$ be fixed is interesting in another
sense.
We obtain non-trivial algebra at $q=0$  if we  normalize
 the generators $T_n$'s
 as $\TT_n=T_n q^{\frac{|n|}{2}}$.
We can define the Verma module generated by $\TT_n$'s with the
highest weight $\lambda$.
It may be seen that at $q=0$ the  Kac determinant
dose not depend on $\lambda$.
Therefore, if $t$ is generic, we have  no singular vectors.
The most remarkable thing at $q=0$ is that after
bosonization every weight
vector can be written down by
the Hall-Littlewood symmetric function\cite{rAKOS}.

\item  The Virasoro algebra plays a central role in the conformal
field theory in two dimensions
 because it generates the conformal transformations of
 the local fields.
 If we want to understand the physical meaning of the
 quantum deformed Virasoro algebra \v which was introduced to study the
$q$-deformed CSM-CFT correspondence,
we have to investigate the symmetry generated by the algebra \v.
In other words, we have to understand  \v from a geometric point of view.
\begin{enumerate}
  \item We have the representation of the  Virasoro algebra for
$c=0$ given by $L_n=-z^{n+1} \partial_z$, i.e. generators of the
conformal transformations of  holomorphic functions.
        It is  interesting to find a similar representation
  for $T_n$'s.
  \item The primary fields in CFT are characterized
by the operator product expansion with the energy
momentum tensor $L(z)$. In the $q$-deformed theory
operator product expansions do not work well. Therefore,
 in order to define  the primary fields of \v,
an algebraic language for treating the intertwining property
 will be needed. In other words we need the coproduct of \v.
  \item The correlation functions of the local fields of the minimal model
satisfy some differential equations\cite{rBPZ}.
What equations do the correlation functions of \v satisfy?
Do that coincide with the $q$-hypergeometric difference equation
\cite{rGR} for some cases?
\item The Jack symmetric functions with  rectangular diagrams
are the singular vectors of  the Virasoro algebra and
diagrams composed of $N$  rectangulars are
 that for  the ${\cal W}_N$ algebra\cite{rAMOS}.
We proved in this  article that  the singular vectors of   \v  are
expressed by
the Macdonald symmetric functions with rectangular diagrams.
Hence we expect that the Macdonald symmetric functions with
diagrams made of $N$ rectangulars correspond to the singular vectors of  a
quantum deformed ${\cal W}_N$ algebra.
    \end{enumerate}
\item A traditional way to introduce the quantum deformation of an infinite
dimensional symmetry algebra such as quantum
affine algebra $U_q(\widehat{\goth{g}})$ \cite{rD} \cite{rJ} is
based on the Yang-Baxter equation.
The principle we have used,
 on the other hand, is just to `` factor'' the Macdonald operator
as $\widehat{D}_{q,t} = \oint \frac{dz}{z}\psi(z)T(z)+\cdots$.
 So far, we have not applied the  idea of the
Yang-Baxter equation directly.
 But, there are many evidences that the algebra  \v has deep connection with
the integrable systems. This point should be clarified.
 \begin{enumerate}
  \item The relation between \v and $U_q(\widehat{\goth{sl}}_2)$
must be clarified.
If we proceed to this direction, a lattice theoretical
interpretation  of \v will be desired.
Can we regard \v as an algebra acting
on the Hilbert space of the $X\!X\!Z$ spin chain model?
\item Frenkel and Reshetikhin defined the $q$-$\cal W$  algebra
 based on  $U_q(\widehat{\goth{g}})$ and pointed out the remarkable
formal resemblance between the bosonized $q$-$\cal W$ currents
 and
the dressed vacuum form (DVF)
in the algebraic Bethe ansatz\cite{rKS}.
 This similar structural coincidence can be observed
in the bosonized \v current $T(z)$ (eq. \eq{b1}).
  \item    The sine-Goldon theory is understood as the
        integrable massive deformation theory of CFT.
Since this theory is integrable, there exist  infinitely many
 conserved charges $\{P_1,P_3, \cdots\}$.
It seems interesting to construct a $q$-deformation
of this scenario. Dose it have to do with the quantum
KdV theory\cite{rBLZ}?
\item In the previous paper\cite{rAMOS} we studied
the generalized matrix model and the
Virasoro constraint.
Construct a matrix model whose constraints  are described by \v.
Is any integrable hierarchy associated with such matrix model?

 \end{enumerate}
\end{itemize}
\vskip 5mm

\noindent{\bf Acknowledgments:}


We would like to thank T.~Eguchi, B.~Feigin, E.~Frenkel,
 T.~Inami, M.~Jimbo, S.~Kato, A.~Kuniba, Y.~Matsuo and T.~Miwa
for discussions and encouragements.
This work is supported in part by Grant-in-Aid for Scientific
Research from Ministry of Science and Culture.


\section*{Appendix}
In this appendix, we show our heuristic derivation of
the bosonized \v current $T(z)$. We start from the
known data obtained from the study of the integral
representation of the Macdonald symmetric functions
\cite{rAMOS}\cite{rAOS}.
That are
\begin{itemize}
\item
The bosonized Macdonald operator:
\be
\widehat{D}_{q,t}=\frac{t^N}{t-1}\oint\frac{dz}{2\pi \i z}
\exp\left\{\ps{n}\frac{1-t^{-n}}{n}a_{-n}z^n\right\}
\exp\left\{-\ps{n}\frac{1-t^{n}}{n}a_{n}z^{-n}\right\}-\frac{1}{t-1}.
\ee
\item The creation operator and the zero-mode parts of the screening currents:
\bea
  S_+(z)
  &\!\!=\!\!&
  \exp\left\{\ps{n}\frac{1-t^n}{1-q^n}\frac{a_{-n}}{n}z^{n}\right\}
  \exp\left\{{\rm annihilation\; part\; for }S_+
  \right\} e^{\beta Q}z^{2\beta a_0},\\
  S_-(z)
  &\!\!=\!\!&
  \exp\left\{-\ps{n}\frac{a_{-n}}{n}z^{n}\right\}
  \exp\left\{{\rm annihilation \;part \;for } S_- \right\}
  e^{- Q}z^{-2 a_0}.
\eea
Note that the zero-mode parts are borrowed from the
$q=1$ free boson realization of the screening currents.
\end{itemize}
We want to ``factorize'' the bosonized Macdonald operator.
First, let us
choose the ansatz for factorization as follows
\begin{eqnarray*}
&&\oint\frac{dz}{2\pi \i z}
\exp\left\{\ps{n}\frac{1-t^{-n}}{n}a_{-n}z^n\right\}
\exp\left\{-\ps{n}\frac{1-t^{n}}{n}a_{n}z^{-n}\right\}\\
&=&\oint\frac{dz}{2\pi \i z}
\exp\left\{\ps{n}\frac{1-t^{-n}}{1+q^{\epsilon n}}\frac{a_{-n}}{n}z^n
q^{\epsilon n}\right\}q^{\delta a_0}
\left(B_+(z)+\eta B_-(z)\right)-\eta q^{(\delta+\gamma)a_0},
\end{eqnarray*}
where
\bea
B_+(z)&=&\exp\left\{\ps{n}\frac{1-t^{-n}}{1+q^{\epsilon
n}}\frac{a_{-n}}{n}z^n\right\}
\exp\left\{-\ps{n}(1-t^{n})\frac{a_{n}}{n}z^{-n}\right\}q^{-\delta a_0},\\
B_-(z)&=&\exp\left\{-\ps{n}\frac{1-t^{-n}}{1+q^{\epsilon n}}\frac{a_{-n}}{n}z^n
q^{\epsilon n}\right\}
\exp\left\{\ps{n}(1-t^{n})\frac{a_{n}}{n}z^{-n} q^{\alpha n}\right\}
q^{\gamma a_0},
\eea
and $\alpha,\delta,\gamma,\eta$ and $\epsilon$ are parameters.
Expand $B_{\pm}(z)$ as $B_{\pm}(z)=\sum B_{\pm,n} z^{-n}$.

The problem is to find suitable expressions for the
annihilation parts of $S_{\pm}(z)$ and fix the parameters
$\alpha,\delta,\gamma,\eta,\epsilon$ by imposing the following
conditions.

{\it
The commutation relations between
$B_{+,n}+\eta B_{-,n}$ and $S_{\pm}(w)$ are given by total differences of
some fields $O_{\pm}(w)$ with some parameters $\chi,\xi$:}
\bea
 \left[B_{+,n}+\eta B_{-,n},S_+(w)\right]&=& \frac{d_{\chi}}{d_{\chi}w}
   O_+(w),
\label{e:cond1}\\
 \left[B_{+,n}+\eta B_{-,n},S_-(w)\right]&=& \frac{d_{\xi}}{d_{\xi}w} O_-(w).
\label{e:cond2}
\eea

The reasons why we applied this ansatz i.e. to factorize the Macdonald
operator
by $B_+(z)+\eta B_-(z)$ are the following. i)
Since the RHS's of \eq{cond1} and \eq{cond2} are the total differences
(that are given by sums of two operators),
we have to introduce ``two terms'' $B_+$ and $B_-$.
This principle worked very well in the study of the
screening currents of the
quantum affine algebra $U_q(\widehat{\goth{sl}}_N)$ \cite{rAOS2}.
ii)
Frenkel and Reshetikhin succeeded in finding the $q$-Virasoro
algebra and its bosonization basing on the
quantum affine algebra $U_q(\widehat{\goth{sl}}_2)$ at the critical
level $k=-2$. Their bosonized formula of the $q$-Virasoro
algebra also consists of two terms as
$S(z)=q^{-1}e^{-\lambda(qz)}+q e^{\lambda(q^{-1}z)}$ (see \cite{rFR}
as for notation).
Therefore, if we start from this ansatz and solve the problem,
there may be a chance to observe close connection with
their theory.
However the factorization problem is quite complicated,
it seems very difficult to find an ansatz having higher symmetry
of the parameters at the beginning as
Frenkel-Reshetikhin's $q$-Virasoro current has.

Let us proceed to find the annihilation parts of $S_{\pm}$.
We have the following operator product expansionsn
\bea
B_+(z)S_+(w)&=&:B_+(z)S_+(w):\frac{z-w}{z-tw}q^{-\delta},\\
B_-(z)S_+(w)&=&:B_-(z)S_+(w):\frac{z-q^{\alpha}tw}{z-q^{\alpha}w}q^{\gamma}.
\eea
So as to obtain \eq{cond1},\eq{cond2}, it is desired that
the commutation relations between $B_{\pm}(z)$ and $S_+(w)$ can be
factorized by the delta functions $\delta(q^{\rho_{\pm}}w/z)$ with some
suitable parameters $\rho_{\pm}$. To this end,
let us impose the following operator product expansions:
\bea
S_+(w)B_+(z)&=&:B_+(z)S_+(w):\frac{w-z}{tw-z}q^{-\delta},\\
S_+(w)B_-(z)&=&:B_-n(z)S_+(w):\frac{q^{\alpha}tw-z}{q^{\alpha}w-z}q^{\gamma}.
\eea
Then we can find the annihilation part of $S_+(z)$ as follows
\be
  S_+(z)
=
  \exp\left\{\ps{n}\frac{1-t^n}{1-q^n}\frac{a_{-n}}{n}z^{n}\right\}
  \exp\left\{-\ps{n}(1+q^{\epsilon n})\frac{1-t^n}{1-q^n}\frac{a_{n}}{n}z^{-n}
  \right\} e^{\beta Q}z^{2\beta a_0},
\ee
and we also have the relations
$\epsilon=-\alpha,\delta=-\beta$ and $\gamma =-\beta$.

Repeat the same argument with $B_{\pm}(z)$ and $S_-(z)$.
We have
\bea
B_+(z)S_-(w)&=&:B_+(z)S_-(w):\frac{z-q w}{z-w}q^{-1},\\
B_-(z)S_-(w)&=&:B_-(z)S_-(w):\frac{z-q^{\alpha}w}{z-q^{\alpha+1}w}q.
\eea
Postulating
\bea
S_-(w)B_+(z)&=&:B_+(z)S_-(w):\frac{q w-z}{w-z}q^{-1},\\
S_-(w)B_-(z)&=&:B_-(z)S_-(w):\frac{q^{\alpha}w-z}{q^{\alpha+1}w-z}q,
\eea
we obtain
\be
  S_-(z)
=
  \exp\left\{-\ps{n}\frac{a_{-n}}{n}z^{n}\right\}
  \exp\left\{\ps{n}(1+q^{\epsilon n})\frac{a_{n}}{n}z^{-n}q^{-n} t^{n}\right\}
  e^{- Q}z^{-2 a_0}.
\ee
Summarizing the commutation relations, we have
\bea
\left[B_+(z),S_+(w)\right]&=&:B_+(z)S_+(w):(t-1)\delta(tw/z),\\
\left[B_-(z),S_+(w)\right]&=&:B_-(z)S_+(w):(t^{-1}-1)\delta(q^{\alpha}w/z),\\
\left[B_+(z),S_-(w)\right]&=&:B_+(z)S_-(w):(q^{-1}-1)\delta(w/z),\\
\left[B_-(z),S_-(w)\right]&=&:B_-(z)S_-(w):(q-1)\delta(q^{\alpha+1}w/z).
\eea
These equations are equivalent to the following:
\bea
&&\left[B_{+,n},S_+(w)\right]=
  \exp\left\{\ps{n}\frac{(1-t^n)(q^{\epsilon n}+q^n)}{(1+q^{\epsilon
n})(1-q^n)}
\frac{a_{-n}}{n}w^{n}\right\} \label{e:ap}\\
&& \times \exp\left\{-\ps{n}\frac{(1-t^n)(t^{-n}+1-q^nt^{-n}+q^{\epsilon
n})}{1-q^n}
\frac{a_{n}}{n}w^{-n}
  \right\} e^{\beta Q}q^{\beta a_0}w^{2\beta a_0}(t-1)t^nw^n,  \n
&&\left[B_{-,n},S_+(w)\right]=\exp\left\{\ps{n}
\frac{(1-t^n)(t^{-n}+1+q^nt^{-n}-q^{\epsilon n})}{(1+q^{\epsilon n})(1-q^n)}
\frac{a_{-n}}{n}w^{n}\right\} \label{e:bp}\\
&& \times  \exp\left\{-\ps{n}\frac{(1-t^n)(q^{\epsilon n}+q^n)}{1-q^n}
\frac{a_{n}}{n}w^{-n}
  \right\} e^{\beta Q}q^{-\beta a_0}w^{2\beta a_0}(t^{-1}-1)p^{-n}w^n,\n
&&\left[B_{+,n},S_-(w)\right]=
\exp\left\{-\ps{n}\frac{t^{-n}+q^{\epsilon n}}{1+q^{\epsilon n}}
\frac{a_{-n}}{n}w^{n}\right\} \label{e:am} \\
&& \times  \exp\left\{\ps{n}(-1+t^n+q^{-n}t^{n}+q^{(\epsilon-1) n}t^n)
\frac{a_{n}}{n}w^{-n}
  \right\} e^{-Q}q^{\beta a_0}w^{-2 a_0}(q^{-1}-1)w^n ,   \n
&& \left[B_{-,n},S_-(w)\right]=
\exp\left\{-\ps{n}\frac{q^n-q^{n}t^{-n}+1+q^{\epsilon n}}{1+q^{\epsilon n}}
\frac{a_{-n}}{n}w^{n}\right\} \label{e:bm} \\
&& \times  \exp\left\{\ps{n}(q^{-n}+q^{(\epsilon-1) n}t^n)
\frac{a_{n}}{n}w^{-n}
  \right\} e^{-Q}q^{-\beta a_0}w^{-2 a_0}(q-1)q^{(\alpha +1)n}w^n.
  \nonumber
\eea

Let us examine in what conditions can we have
\eq{cond1},\eq{cond2}.
Studying the zero-mode factor in the equations \eq{ap} and \eq{bp},
namely $q^{\beta a_0}w^{2\beta a_0}$ and $q^{-\beta a_0}w^{2\beta
  a_0}$, it is found that we have to set $\chi=q$. Next, examining the
oscillator
factors, we have
\bea
&&\frac{(1-t^n)(q^{\epsilon n}+q^n)}{(1+q^{\epsilon n})(1-q^n)}=
\frac{(1-t^n)(t^{-n}+1+q^nt^{-n}-q^{\epsilon n})}{(1+q^{\epsilon
n})(1-q^n)}q^n\\
&&\frac{(1-t^n)(t^{-n}+1-q^nt^{-n}+q^{\epsilon n})}{1-q^n}=
\frac{(1-t^n)(q^{\epsilon n}+q^n)}{1-q^n}q^{-n}
\eea
These equations can be solved uniquely and we obtain
$q^{\epsilon}=qt^{-1}
\equiv p$.
In the same way,
we have $\xi=t$ from the zero-mode factors in the
equations
\eq{am},\eq{bm} and condition from the oscillator factors
gives us also $q^{\epsilon}=p$.
We have obtained
\be
q^{-\alpha}=q^{\epsilon}=p,\;\;q^{\delta}=q^{\gamma}=t^{-1},\;\;\;
\chi=q,\;\;\xi=t.
\ee

Finally if we set $\eta=p^{-1}$ then
we have the desired equations
\bea
 \left[B_{+,n}+p^{-1} B_{-,n},S_+(w)\right]&=& \frac{d_{q}}{d_{q}w}
   O_+(w),\\
 \left[B_{+,n}+p^{-1} B_{-,n},S_-(w)\right]&=& \frac{d_{t}}{d_{t}w} O_-(w),
\eea
where
\bea
  O_+(w)
  &\!\!=\!\!&-(1-q)(1-t^{-1})(p^{-1} w)^{n+1}
  \exp\left\{\ps{n}\frac{1+t^n}{1+p^n}\frac{1-t^n}{1-q^n}
  \frac{a_{-n}}{n} w^nt^{-n}\right\} \\&&\times
  \exp\left\{-\ps{n}(1+t^n)\frac{1-t^n}{1-q^n}\frac{a_{n}}{n}
  w^{-n}p^n\right\}
  e^{\beta Q }w^{2\beta a_0}t^{-a_0},
 \\
  O_-(w)
  &\!\!=\!\!& -(1-q^{-1})(1-t)w^{n+1}
  \exp\left\{-\ps{n} \frac{1+q^n}{1+p^n}\frac{a_{-n}}{n}
  w^nt^{-n}\right\} \\&&\times
  \exp\left\{\ps{n} (1+q^n)\frac{a_{n}}{n}
  w^{-n}p^{-n}\right\}e^{-Q}w^{-2 a_0}t^{a_0}.
\eea

%
%


\begin{thebibliography}{99}
%
\bibitem{rCS}F.~Calogero, {\sl Jour. Math. Phys.}
{\bf 10} (1969) 2197-2200;\\
B.~Sutherland, {\sl Jour. Math. Phys.} {\bf 12}  (1970) 246-250, 251-256.
%
\bibitem{rBPZ}
  A.A.~Belavin, A.M.~Polyakov and A.B.~Zamolodchikov,
 {\sl Nucl. Phys.} {\bf B241} (1984) 333-380.
%
\bibitem{rAMOS}
  H.~Awata, Y.~Matsuo, S.~Odake and J.~Shiraishi,
  {\it Collective Field Theory, Calogero-Sutherland Model
   and Generalized Matrix Models},
  {\sl Phys. Lett.} {\bf B347} (1995) 49-55.
  {\it A Note on Calogero-Sutherland Model, $W_n$ Singular Vectors
   and Generalized Matrix Models},
  preprint, hep-th/9503028,
  to appear in  Soryushiron kenkyu (Kyoto);
  {\it Excited State of Calogero-Sutherland Model and Singular
Vectors of the  $W_n$ Algebra},  preprint, hep-th/9503043,
to appear in {\sl Nucl. Phys. {\bf B}}.
%
\bibitem{rAOS}
 H.~Awata, S.~Odake and J.~Shiraishi,
{\it Integral Representations of the Macdonald Symmetric Functions
   and Generalized Matrix Models},
  preprint, q-alg/9506006
%
\bibitem{rMY1}
  K.~Mimachi and Y.~Yamada,
  {\it Singular vectors of the Virasoro algebra in terms of
  Jack symmetric polynomials},
  preprint (November 1994).
%
\bibitem{rS} R.~Stanley, {\sl Adv. Math.} {\bf 77}
(1989) 76-115;\\
I.G.~Macdonald, Lect. Note in Math. {\bf 1271},
189-200, Springer 1987.
%
\bibitem{rFF}
  B.~L.~Feigin and D.~B.~Fuchs, {\it Skew-symmetric
differential operators on the line and Verma modules
 over the Virasoro algebra}, {\sl Functs. Anal. Prilozhen}
{\bf16} (1982) 47;
 {\it Verma modules over the Virasoro algebra},
in L.~D.~Faddeev nad A.~A.~Malcev, eds.
 {\sl Topology, Proceedings of Leningrad conference}, 1982,
 Lect. Note in Math.
 {\bf 1060}, Springer 1985.
%
\bibitem{rMac} I.G.~Macdonald
``Symmetric Functions and Hall Functions'',
Oxford University Press 1979.
%
\bibitem{rJS}
  A.~Jevicki and B.~Sakita,
  {\sl Nucl. Phys.} {\bf B165} (1980) 511-527.
%
\bibitem{rAJL} I.~Andri\'c, A.~Jevicki and H.~Levin,
{\sl Nucl. Phys.} {\bf B215 [FS7]} (1983) 307-315.
%
\bibitem{rDFF}
  G.~Felder,
  {\sl Nucl. Phys.} {\bf B317} (1989) 215-236,
  Errata {\sl Nucl. Phys.} {\bf B324} (1989) 548.
%
\bibitem{rFR}E.~Frenkel and N.~Reshetikhin,
{\em Quantum Affine Algebras and Deformations
of The Virasoro and $\cal W$-Algebra}, q-alg/9505025 May 1995.
%
\bibitem{rAKOS}H.~Awata, H.~Kubo, S.~Odake and J.~Shiraishi,
 {\it in preparation}.
%
\bibitem{rD}V.~G.~Drinfeld,
{\it Quantum Groups}, {\sl ICM} {\bf 86} {\sl report}.
%
\bibitem{rJ}M.~Jimbo, {\it A $q$-Difference Analogue of
$U(\goth{g})$ and
 the Yang-Baxter Equation}, {\sl Lett. Math. Phys.} {\bf 10} (1985) 63-69.
%
\bibitem{rGR}
  G.~Gasper and M.~Rahman,
  {\it Basic Hypergeometric Series},
  Encyclopedia of Mathematics and its Applications,
 Cambridge University Press 1990.

\bibitem{rKS}
  A.~Kuniba and J.~Suzuki,
  {\it Analytic Bethe Ansatz for Fundamental Representations of
Yangians}, preprint hep-th/9406180.
%
\bibitem{rBLZ}
V.~Bazhanov, S.~Lukyanov, A.~Zamolodchikov,
{\it Integrable Structure of Conformal Field Theory, Quantum KdV Theory and
Thermodynamic Bethe Ansatz}, preprint hep-th/9412229.
%

\bibitem{rAOS2}
  J.~Shiraishi,
  {\it Free Boson Representation of $U_q(\widehat{\goth{sl}}_2)$},
  {\sl Phys. Lett.} A {\bf 171}, (1992) 243-248,\\
  H.~Awata, S.~Odake and J.~Shiraishi,
  {\it Free Boson Representation of $U_q(\widehat{\goth{sl}}_3)$},
  {\sl Lett. Math. Phys.} {\bf 30}, (1994) 207-216;
  {\it Free Boson Realization of $U_q(\widehat{\goth{sl}}_N)$},
  {\sl Commun. Math. Phys.} {\bf 162}, (1994) 61-83.



\end{thebibliography}
\end{document}